\documentclass[12pt]{article}
\usepackage{amsmath}
\usepackage{amssymb}
\usepackage{amsthm}
\usepackage{fullpage}
\usepackage[utf8]{inputenc}
\usepackage{mathtools}
\usepackage{url}
\usepackage[usenames,svgnames]{xcolor}

\usepackage[pagebackref]{hyperref} 
\usepackage[all]{hypcap} 

\theoremstyle{plain}
\newtheorem{theorem}{Theorem}[section]

\theoremstyle{definition}

\newtheorem{remark}[theorem]{Remark}

\numberwithin{equation}{section} 

\hypersetup{
  breaklinks,colorlinks,
  citecolor=Green,
  linkcolor=BlueViolet,
  urlcolor=DarkBlue,
  pdftitle={On  some class of partial difference equations admitting a zero-curvature representation},
  pdfauthor={Svinin Andrei K.},
}

\begin{document}
\baselineskip 16pt

\medskip
\begin{center}
\begin{Large}\fontfamily{cmss}
\fontsize{17pt}{27pt}
\selectfont
\textrm{Tuenter polynomials and a Catalan triangle}
\end{Large}\\
\bigskip
\begin{large} {Andrei K. Svinin}
 \end{large}
\\
\bigskip
\begin{small}
{\em Matrosov Institute for System Dynamics and Control Theory of 
Siberian Branch of Russian Academy of Sciences}\\
svinin@icc.ru \\
\end{small}
\end{center}
\bigskip

\begin{abstract}
We consider Tuenter polynomials as linear combinations of  descending factorials and  
show that  coefficients of these linear combinations are expressed via a Catalan triangle of numbers. We also describe a triangle of coefficients in terms of some polynomials. 
\end{abstract}


\section{Preliminaries. Tuenter polynomials}
The polynomials we are going to study in this brief note are defined by a recursion \cite{Tuenter}
\begin{equation}
P_{k+1}(n)=n^2\left(P_k(n)-P_{k}(n-1)\right)+nP_k(n-1),\;\;
n\in\mathbb{N}
\label{1}
\end{equation}
with initial condition $P_0(n)=1$. The first few polynomials yielded by (\ref{1}) are as follows.
\[
P_1(n)=n,
\]
\[
P_2(n)=n(2n-1),
\]
\[
P_3(n)=n(6n^2-8n+3),
\]
\[
P_4(n)=n(24n^3-60n^2+54n-17),
\]
\[
P_5(n)=n(120n^4-480n^3+762n^2-556n+155),
\]
\[
P_6(n)= n(720n^5-4200n^4+10248n^3-12840n^2+8146n-2073).
\]
Let us refer to these polynomials as Tuenter ones. Introducing a recursion operator $R:=n^2\left(1-\Lambda^{-1}\right)+n\Lambda^{-1}$, where $\Lambda$ is a shift operator acting as $\Lambda(f(n))=f(n+1)$, one can write $P_k(n)=R^k(1)$.
The sense of these polynomials is that they help to count the sum
\[
S_{r}(n)=\sum_{j=0}^{2n}\left(
\begin{array}{c}
2n\\
j
\end{array}
\right)|n-j|^r
\]
for odd $r$. 

Bruckman in \cite{Bruckman} asked to prove that $S_3(n)=n^2\left(
\begin{array}{c}
2n\\
n
\end{array}
\right)$. Strazdins in \cite{Strazdins} solved this problem and conjectured that   $S_{2k+1}(n)=\tilde{P}_{k}(n)\left(
\begin{array}{c}
2n\\
n
\end{array}
\right)$ with some monic polynomial $\tilde{P}_{k}(n)$ for any $k\geq 0$. Tuenter showed in \cite{Tuenter} that it is almost true. More exactly, he proved that
\[
S_{2k+1}(n)=P_{k}(n)n\left(
\begin{array}{c}
2n\\
n
\end{array}
\right)=P_{k}(n)\frac{(2n)!}{(n-1)!n!}.
\]
One can see that polynomial $\tilde{P}_{k}(n)$ is monic only for $k=0, 1$. The recursion (\ref{1}) follows from  \cite{Tuenter}
\[
S_{r+2}(n)=n^2S_r(n)-2n(2n-1)S_r(n-1).
\]
Also, as was noticed in \cite{Tuenter}, polynomials $P_{k}(n)$ can be obtained as a special case of Dumont-Foata polynomials of three variables \cite{Dumont}.

\section{The Tuenter polynomials as linear combinations of descending factorials }
Consider descending factorials 
\[
(n)_k:=n(n-1)(n-2)\cdots(n-k+1).
\]
It can be easily seen that
\begin{equation}
R((n)_k)=k^2(n)_k+(k+1)(n)_{k+1}.
\label{2}
\end{equation}
Let us consider $P_k(n)$ as linear combinations of descending factorials
\[
P_k(n)=\sum_{j=1}^kc_{j,k}(n)_j,
\]
with some coefficients $c_{j,k}$  to be calculated. For example, for the first few $P_k(n)$ we get
\[
P_1(n)=(n)_1,\;\;
\]
\[
P_2(n)=(n)_1+2(n)_2,\;\;
\]
\[
P_3(n)=(n)_1+10(n)_2+6(n)_3,
\]
\[
P_4(n)=(n)_1+42(n)_2+84(n)_3+24(n)_4,
\]
\[
P_5(n)=(n)_1+170(n)_2+882(n)_3+720(n)_4+120(n)_5,
\]
\[
P_6(n)=(n)_1+682(n)_2+8448(n)_3+15048(n)_4+6600(n)_5+720(n)_6.
\]
With (\ref{2}) we can easily derive recurrence relations for the coefficients $c_{j,k}$. Indeed, from
\begin{eqnarray}
P_{k+1}(n)&=&\sum_{j=1}^{k+1}c_{j,k+1}(n)_j\nonumber\\
&=&R(P_k(n))\nonumber\\
&=&\sum_{j=1}^kc_{j,k}\left(j^2(n)_j+(j+1)(n)_{j+1}\right)\nonumber
\end{eqnarray}
we get 
\begin{equation}
c_{j,k+1}=j^2c_{j, k}+jc_{j-1,k},\;\;j\geq 1,\;\;k\geq j.
\label{3}
\end{equation}
To use (\ref{3}), one must  agree that  $c_{0,k}=c_{k+1,k}=0$ for $k\geq 1$. Then, starting from $c_{1,1}=1$ we obtain  the whole set $\{c_{j,k} : j\geq 1, k\geq j\}$. For example, $c_{1,k}=1$ for all $k\geq 1$, while for $j=2$ we obviously get a recursion
\[
c_{2,k+1}=4c_{2,k}+2,\;\;c_{2,1}=0.
\]
As can be easily seen, a solution of this equation is given by 
\begin{equation}
c_{2,k}=\frac{1}{3}\left(2^{2k-1}-2\right),\;\;
k\geq 2.
\label{5}
\end{equation}
\begin{remark}
It is interesting to note that  integer sequence (\ref{5}), known as \textrm{A020988} in \cite{Sloane} gives $n$-values of local maxima for $s(n):=\sum_{j=1}^na(j)$, where $\{a(n)\}$ is the Golay-Rudin-Shapiro sequence \cite{Brillhart}. 
\end{remark}

For  the whole set of the coefficients $\left\{c_{j,k}\right\}$, we get the following.
\begin{theorem} \label{th1}
A solution of equation (\ref{3}) with $c_{0,k}=c_{k+1,k}=0$ for $k\geq 1$ and $c_{1,1}=1$ is given by
\begin{equation}
c_{j,k}=\frac{j!}{(2j-1)!}\left(\sum_{q=1}^j (-1)^{q+j}B_{j,q}q^{2k-1} \right),\;\;\forall j\geq 1,\;\;
k\geq j,
\label{4}
\end{equation}
where the numbers
\[
B_{j,q}:=\frac{q}{j}\left(
\begin{array}{c}
2j \\
j-q
\end{array}
\right)
\]
constitute a Catalan triangle \cite{Shapiro}.
\end{theorem}

\noindent
\textbf{Proof}.
Substituting (\ref{4}) into (\ref{3}) and collecting terms at $q^{2k-1}$, we obtain that sufficient condition for (\ref{4}) to be a solution of (\ref{3}) is that the numbers $B_{j,q}$  enjoy the relation
\[
\frac{q^2 j!}{(2j-1)!}B_{j,q}=\frac{j^2 j!}{(2j-1)!}B_{j,q}-\frac{j!}{(2j-3)!}B_{j-1,q},\;\;\forall q=1,\ldots, j-1.
\]
Simplifying the latter we get the relation
\[
(j-q)(j+q)B_{j,q}=(2j-1)(2j-2)B_{j-1,q}
\]
which can be easily verified. Therefore the theorem is proved. $\Box$

The set $\{c_{j, k}\}$ can be presented as the number triangle  
\[
\begin{array}{cccccc}
          &          &c_{1, 1} &            &   & \\
          &c_{1, 2}  &         &c_{2, 2}    &    & \\
c_{1, 3}  &          &c_{2, 3} &            &c_{3, 3}& \\
  &   \ddots       & & \ddots           & &\ddots
\end{array}
\]
whose description is given by theorem \ref{th1}.

\begin{remark}
From \cite{Shapiro} one knows that the number $B_{j,q}$ can be interpreted as the number of pairs of non-intersecting  paths of length $j$ and distance $q$. The Catalan numbers itself (\textrm{A000108})  are 
\[
C_j:=B_{j,1}=\frac{1}{j}\left(
\begin{array}{c}
2j \\
j-1
\end{array}
\right).
\]
\end{remark}

Therefore, we got an infinite number of integer sequences each of which is defined by numbers from the Catalan triangle and  begins from $c_{j,j}=j!$. 
Let us list the first few ones. For example, one has, 
\[
c_{1, k}=1,
\]
\[
c_{2,k}=\frac{1}{3}\left(2^{2k-1}-2\right),
\]
\[
c_{3,k}=\frac{1}{20}\left(3^{2k-1}-4\cdot 2^{2k-1}+5\right),
\]
\[
c_{4,k}=\frac{1}{210}\left(4^{2k-1}-6\cdot 3^{2k-1}+14\cdot 2^{2k-1}-14\right),
\]
\[
c_{5,k}=\frac{1}{3024}\left(5^{2k-1}-8\cdot 4^{2k-1}+27\cdot 3^{2k-1}-48\cdot 2^{2k-1}+42\right),
\]
\[
c_{6,k}=\frac{1}{55440}\left(6^{2k-1}-10\cdot 5^{2k-1}+44\cdot 4^{2k-1}-110\cdot 3^{2k-1}+165\cdot 2^{2k-1}-132\right),\ldots
\]
All these sequences are indeed integer because they are solutions of (\ref{3}). 

 Let us replace $k\mapsto j+k$ in (\ref{3}) and seek its solution in the form
$c_{j, j+k}=F_k(j)j!$. Substituting the latter in (\ref{3}) we come to the recurrence relation
\begin{equation}
F_k(j)-F_{k}(j-1)=j^2F_{k-1}(j)
\label{rec1}
\end{equation}
with conditions $F_0(j)=1$ and $F_k(1)=1$. A solution of (\ref{rec1}) is
\[
F_k(j)=1+\sum_{2\leq\lambda_1\leq j}\lambda_1^2+\sum_{2\leq\lambda_1\leq\lambda_2\leq j}\lambda_1^2\lambda_2^2+\cdots+\sum_{2\leq\lambda_1\leq\cdots\leq \lambda_k\leq j}\lambda_1^2\lambda_2^2\cdots\lambda_k^2.
\]
In particular,
\[
F_1(j)=1+2^2+\cdots+j^2=\frac{1}{6}j(j+1)(2j+1),
\]
that is, $F_1(j)$ yields \textrm{A000330} sequence of the square pyramidal numbers. The next two  polynomials $F_k(j)$ are
\[
F_2(j)=\frac{1}{360}j(j+1)(j+2)(2j+1)(2j+3)(5j-1)
\]
and
\[
F_3(j)=\frac{1}{45360}j(j+1)(j+2)(j+3)(2j+1)(2j+3)(2j+5)(35j^2-21j+4).
\]
Looking at these examples and others we can suppose that 
\[
F_k(j)=\prod_{q=0}^k(j+q)\prod_{q=0}^{k-1}(2j+2q+1)\tilde{F}_k(j),
\]
where $\tilde{F}_k(j)$ is a polynomial of $(k-1)$ degree, which satisfy 
\[
\tilde{F}_k(1)=\frac{1}{(k+1)!(2k+1)!!}.
\]

\section*{Acknowlegments}
This  work  was  supported  in  part  by the Council for Grants of the
 President  of  Russian  Foundation  for  state  support of the leading
 scientific schools, project NSh-8081.2016.9.


\begin{thebibliography}{10}

\bibitem{Brillhart}
J. Brillhart and P. Morton, A case study in mathematical research: the Golay-Rudin-Shapiro sequence, The American Mathematical Monthly, \textbf{103} (1996), 854-869. 

\bibitem{Bruckman}
P.~S. Bruckman, Problem B-871, Fibonacci Quartely, \textbf{37} (1999), 85.

\bibitem{Dumont}
D. Dumont and D. Foata, Une propri\'et\'e de sym\'etrie des nombres de Genocchi, Bulletin de la Soci\'et\'e Math\'ematique de France \textbf{104} (1976), 433-451. 

\bibitem{Shapiro}
L.~W. Shapiro, A Catalan triangle, Discrete Math., \textbf{14} (1976), 83-90.  

\bibitem{Sloane}
N.~J.~A. Sloane. The On-Line Encyclopedia of Integer Sequences, http://oeis.org.

\bibitem{Strazdins}
I. Strazdins, Solution to problem B-871, Fibonacci Quartely, \textbf{38.1} (2000), 86-87.


\bibitem{Tuenter}
H.~J.~H. Tuenter, Walking into an absolute sum, Fibonacci Quartely, \textbf{40} (2002), 175-180. 
 
\end{thebibliography}

\end{document}